# An introduction to the theory of citing[*]

M. V. Simkin and V. P. Roychowdhury

*Department of Electrical Engineering, University of California, Los Angeles, CA 90095-1594*

## Read before you cite

Many psychological tests have the so-called *lie-scale*. A small but sufficient number of questions that admit only one true answer, such as "Do you **always** reply to letters immediately after reading them?" are inserted among others that are central to the particular test. A wrong reply for such a question adds a point on the lie-scale, and when the *lie-score* is high, the over-all test results are discarded as unreliable. Perhaps, for a scientist the best candidate for such a lie-scale is the question "Do you read **all** of the papers that you cite?"

Comparative studies of the popularity of scientific papers have been a subject of much interest [1]-[4], but the scope was limited to citation counting. We discovered a method of estimating what percentage of people who cited the paper have actually **read** it. Remarkably, this can be achieved without any testing of the scientists, but solely on the basis of the information available in the ISI citation database.

Freud [5] discovered that the application of his technique of Psychoanalysis to slips in speech and writing could reveal a lot of hidden information about human psychology. Similarly, we find that the application of statistical analysis to **misprints in scientific citations** can give an insight into the process of scientific writing. As in the Freudian case, the truth revealed is embarrassing. For example, an interesting statistic revealed in our study is that a lot of misprints are identical. The probability of repeating someone else's misprint accidentally is small. One concludes that repeat misprints are most likely due to copying from a reference list used in another paper.

Our initial report [6] led to a lively discussion [7] on whether copying citations is tantamount to not reading the original paper. Alternative explanations are worth exploring; however, such hypotheses should be supported by data and not by anecdotal claims. It is indeed most natural to assume that a copying citer also failed to read the paper in question (albeit this can not be rigorously proved). *Entities must not be multiplied beyond necessity.* Having thus shaved the critique with Occam's razor, we will proceed to use the term non-reader to describe a citer who copies.

As misprints in citations are not too frequent, only celebrated papers provide enough statistics to work with. Let us have a look at the distribution of misprints in citations to one renowned paper [8], which accumulated 4300 citations. Out of these citations 196 contain misprints, out of which only 45 are distinct. The most popular misprint in a page number appeared 78 times.

As a preliminary attempt, one can estimate the ratio of the number of readers to the number of citers, $R$, as the ratio of the number of **distinct** misprints, $D$, to the **total number** of misprints, $T$. Clearly, among $T$ citers, $T - D$ copied, because they repeated someone else's misprint. For the $D$ others, with the information at hand, we have no evidence that they did not read, so according to the presumed innocent principle, we assume that they did. Then in our sample, we have $D$ readers and $T$ citers, which lead to:

$$R \approx D/T. \qquad (1)$$

Substituting $D = 45$ and $T = 196$ in Eq.(1), we obtain that $R \approx 0.23$. This estimate would be correct if the people who introduced original misprints had always read the original paper. It is more reasonable to assume that the probability of introducing a new misprint in a citation does not depend on whether the author had read the original paper. Then, if the fraction of read citations is $R$, the number of readers in our sample is $RD$, and the ratio of the number of readers to the number of citers in the sample is $RD/T$. What happens to our estimate, Eq. (1)? It is correct, just the sample is not representative: the fraction of read citations among the citations containing misprints is less than in the general citation population.

Can we still determine $R$ from our data? Yes. From the misprint statistics we can determine the average number of times, $n_p$, a typical misprint propagates:

$$n_p = \frac{T - D}{D}. \qquad (2)$$

The number of times a misprint had propagated is the number of times the citation was copied from either the paper which introduced the original misprint, or from one of subsequent papers, which copied (or copied from copied etc) from it. A misprinted citation should be no different from a correct citation as far as copying is concerned. This means that a selected at random citation, on average, is copied (including copied from copied etc) $n_p$ times. The read citations are no different from unread citations as far as copying goes. Therefore, every read citation, on average, was copied $n_p$ times. The fraction of read citations s is thus

$$R = \frac{1}{1+n_p}. \quad (3)$$

After substituting Eq. (2) into Eq. (3), we recover Eq. (1).

Note, however, that the average number of times a misprint propagates is not equal to the number of times the citation was copied, but to the number of times it was copied *correctly*. Let us denote the average number of citations copied (including copied from copied etc) from a particular citation as $n_c$. It can be determined from $n_p$ the following way. The $n_c$ consists of two parts: $n_p$ (the correctly copied citations) and misprinted citations. If the probability of making a misprint is $M$ and the number of correctly copied citations is $n_p$ then the total number of copied citations is $\frac{n_p}{1-M}$ and the number of misprinted citations is $\frac{n_p M}{1-M}$. As each misprinted citation was itself copied $n_c$ times, we have the following self-consistency equation for $n_c$:

$$n_c = n_p + n_p \times \frac{M}{1-M} \times (1+n_c) \quad (4)$$

Eq. (4) has the solution

$$n_c = \frac{n_p}{1-M-n_p \times M} \quad (5)$$

After substituting Eq. (2) into Eq. (5) we get:

$$n_c = \frac{T-D}{D-MT}. \quad (6)$$

From this we get:

$$R = \frac{1}{1+n_c} = \frac{D}{T} \times \frac{1-(MT)/D}{1-M} \quad (7)$$

The probability of making a misprint can be estimated as $M = \frac{D}{N}$, where $N$ is the total number of citations. After substituting this into Eq. (7) we get:

$$R = \frac{D}{T} \times \frac{N-T}{N-D}. \quad (8)$$

Substituting $D=45$, $T=196$, and $N=4300$ in Equation (8), we get $R \approx 0.22$, which is very close to the initial estimate, obtained using Eq.(1).

## Copied citations create renowned papers

During the "Manhattan project" (the making of nuclear bomb), Fermi asked Gen. Groves, the head of the project, what is the definition of a "great" general [9]. Groves replied that any general who had won five battles in a row might safely be called great. Fermi then asked how many generals are great. Groves said about three out of every hundred. Fermi conjectured that considering that opposing forces for most battles are roughly equal in strength, the chance of winning one battle is ½ and the chance of winning five battles in a row is $1/2^5 = 1/32$. "So you are right General, about three out of every hundred. Mathematical probability, not genius." The existence of military genius was also questioned on basic philosophical grounds by Tolstoy [10].

A commonly accepted measure of "greatness" for scientists is the number of citations to their papers [2]. For example, SPIRES, the High-Energy Physics literature database, divides papers into six categories according to the number of citations they receive. The top category, "Renowned papers" are those with 500 or more citations. Let us have a look at the citations to roughly 24 thousands papers, published in Physical Review D in 1975-1994 [11]. As of 1997 there where about 350 thousands of such citations: fifteen per published paper on the average. However, forty-four papers were cited five hundred times or more. Could this happen if **all papers are created equal**? If they indeed are then the chance to win a citation is one in 24,000. What is the chance to win 500 cites out of 350,000? The calculation is slightly more complex than in the militaristic case, but the answer is one in $10^{500}$, or, in other words, it is zero. One is tempted to conclude that those forty-four papers, which achieved the impossible, are great.

A more careful analysis puts this conclusion in doubt. We just have shown that the majority of scientific citations are copied from the lists of references used in another papers. This way a paper that already was cited is likely to be cited again, and after it is cited again it is even more likely to be cited in the future. In other words, "unto every one that hath shall be given, and he shall have abundance" [12]. This phenomenon is known as "Matthew effect" [13], "cumulative advantage" [14], or "preferential attachment" [15].

The effect of citation copying on the probability distribution of citations can be quantitatively understood within the framework of **the model of random-citing scientists (RCS)** [16], which is as follows. When a scientist is writing a manuscript he picks up $m$ random articles, cites them, and also copies some of their references, each with probability $p$.

This model was stimulated by the recursive literature search model [17] and can be solved using methods developed to deal with multiplicative stochastic processes [18]. These methods are too complicated to be described in a popular article so we will just state the results. A good agreement between the RCS model and actual citation data [7] is achieved with the input parameters $m=3$ and $p=1/4$ (see Figure 1). Now what is the probability for an arbitrary paper to become "renowned", i.e. receive more than five hundred citations? A calculation shows that this probability is one in 600. This means that about 40 out of 24,000 papers should be renowned; ergo, mathematical probability, not genius.

In one incident [19] Napoleon (incidentally, he was the military commander, whose genius was questioned in Ref. [10]) said to Laplace "They tell me you have written this large book on the system of the universe, and have never even mentioned its Creator." The reply was "I have no need for this hypothesis". It is worthwhile to note that Laplace was not against God. He simply did not need to postulate his existence in order to explain existing astronomical data. Similarly, the present work is not blasphemy. Of course, in some spiritual sense, great scientists do exist. It is just that even if they would not exist, citation data would look the same.

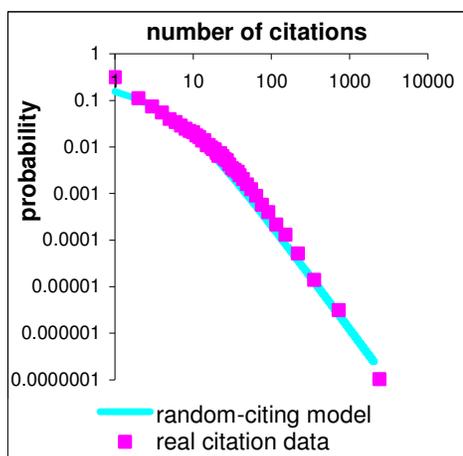

**Figure 1.** The outcome of the model of random-citing compared to actual citation data. Mathematical probability rather than genius can explain why some papers are cited a lot more than the others.

---

* This article (under the different title of "Do you sincerely want to be cited? Or: read before you cite") appeared in the December 2006 issue of *Significance*, a scientific-popular journal published by Royal Statistical Society.